\documentclass[11pt]{article}

\usepackage[utf8]{inputenc}
\usepackage{amsmath}
\usepackage{amsfonts,amssymb,latexsym,multicol,xcolor}
\usepackage[francais]{babel}
\usepackage[T1]{fontenc}

\newcounter{fig}
\newtheorem{theo}{Th\'eor\`eme}

\newtheorem{prop}{Proposition}

\newcommand{\virg}{\raisebox{.7mm}{,}}

\newcommand{\expli}[1]{\;\text{\footnotesize (#1)}}

\makeatletter
\ifcase\@ptsize

\or

\or

\fi
\makeatother

\newcommand{\eps}{\varepsilon}
\newcommand{\fhi}{\varphi}
\newcommand{\ioe}{\leqslant}
\newcommand{\soe}{\geqslant}
\renewcommand{\le}{\leqslant}
\renewcommand{\ge}{\geqslant}
\newcommand{\vers}{\rightarrow}

\newcommand{\demi}{{\frac{1}{2}}}

\newcommand{\Scal}{{\mathcal S}}

\newcommand{\Wcal}{{\mathcal W}}

\newcommand{\Fgot}{{\mathfrak F}}

\newcommand{\cgot}{{\mathfrak c}}

\newcommand{\norm}[1]{\left\| #1 \right\|}

\newcommand{\Nat}{{\mathbb N}}
\newcommand{\Int}{{\mathbb Z}}
\newcommand{\Rat}{{\mathbb Q}}
\newcommand{\Real}{{\mathbb R}}
\newcommand{\Com}{{\mathbb C}}

\newcommand{\fin}{\hfill$\Box$}
\newcommand{\dem}{\noindent {\bf D\'emonstration\ }}
\newcommand{\fine}{\tag*{\mbox{$\Box$}}}

\providecommand{\bysame}{\leavevmode ---\ }
\providecommand{\og}{``}
\providecommand{\fg}{''}
\providecommand{\smfandname}{et}

\newcommand{\W}{\mathcal{W}}

\title{Sur certaines \'equations fonctionnelles approch\'ees,\\ li\'ees \`a la transformation de Gauss}

\author{Michel Balazard et Bruno Martin\thanks{Bruno Martin is supported by ANR Grant MUDERA ANR-14-CE34-0009.}}

\begin{document}
\maketitle

\begin{center}
  {\sc Abstract}
\end{center}
\begin{quote}
{\footnotesize In the line of classical work by Hardy, Littlewood and Wilton, we study a class of functional equations involving the Gauss transformation from the theory of continued fractions. This allows us to reprove, among others, a convergence criterion for a diophantine series considered by Chowla, and to give additional information about the sum of this series.
}
\end{quote}

\begin{center}
  {\sc Keywords}
\end{center}
\begin{quote}
{\footnotesize Fractional part, Continued fractions, Gauss transformation, Approximate functional equations \\MSC classification : 40A05 (11A55)}
\end{quote}



\section{Introduction}

D\'esignons par $\{x\}$ la partie fractionnaire du nombre r\'eel $x$, et par $\lfloor x\rfloor$ sa partie enti\`ere, de sorte que
\[
x=\lfloor x\rfloor + \{x\} \quad ; \quad \lfloor x\rfloor \in \Int \quad ; \quad  0 \ioe \{x\} <1.
\]
\`A partir de la partie fractionnaire, nous d\'efinissons \'egalement la \emph{premi\`ere fonction de Bernoulli normalis\'ee},
\begin{equation*}
 B_1(x)=
 \begin{cases}
  \{x\}-1/2  & (x\notin \Int)\\
  \qquad\! 0 & (x\in \Int),
  \end{cases}
\end{equation*}
et la \emph{transformation de Gauss}, 
\[
\alpha(x)=\{1/x\} \quad (0<x < 1).
\]
Les fonctions it\'er\'ees de $\alpha$ sont d\'efinies en posant 
\begin{align*}
\alpha_0(x)&=x\\
\alpha_k(x)&=\alpha \big(\alpha_{k-1}(x)\big) \quad (k\ge 1).
\end{align*}
Lorsque $x$ est rationnel, $\alpha_k(x)$ est d\'efini tant que $k\ioe K$, o\`u $K$ est un entier, appel\'e \emph{profondeur} de~$x$, \'egal \`a la longueur de sa fraction continue (la d\'efinition pr\'ecise est donn\'ee au \S\ref{profondeur_rationnel} ci-dessous) ; on a alors $\alpha_K(x)=0$.

Nous introduisons \'egalement  les fonctions 
\begin{equation*}\label{t63}
 \beta_k(x)=\alpha_0(x)\alpha_1(x)\cdots\alpha_{k}(x)
\end{equation*}
(par convention, $\beta_{-1}=1$). Si $x$ est rationnel, de profondeur $K$, $\beta_k(x)$ est donc d\'efini si $k\ioe K$.

\smallskip

L'\'etude pr\'esent\'ee dans cet article a pour origine la question de la convergence des s\'eries
\begin{align} 
\fhi_1(x)&= \sum_{n\ge 1} \frac{ B_1(nx)}{n}\label{t208}\\
\Wcal(x)&=\sum_{k\soe 0}(-1)^k\beta_{k-1}(x)\log\big(1/\alpha_k(x)\big),\label{t209}
\end{align}
o\`u $x\in \Real$. Pour la s\'erie $\Wcal(x)$, nous adoptons les conventions suivantes : on pose $\W(0)=~0$,~$\Wcal(x)$ est limit\'ee \`a la somme des termes d'indices $<K$ si~$x$ est un nombre rationnel de profondeur~$K$ appartenant \`a l'intervalle $]0,1[$ (par exemple, $\W(1/k)=\log k$ si $k$ est entier, $k\soe 2$),  et $\Wcal$ est prolong\'ee \`a $\Real$ par p\'eriodicit\'e : $\Wcal(x)=\Wcal(\{x\})$. 

\smallskip

Ces deux s\'eries ont \'et\'e consid\'er\'ees dans les ann\'ees trente du si\`ecle dernier. La premi\`ere apparaît dans l'article de Chowla \cite{chowla}, p. 546, mais elle est tr\`es voisine d'un exemple donn\'e par Riemann au \S XIII de son m\'emoire d'habilitation sur les s\'eries trigonom\'etriques. La seconde apparaît dans l'article de Wilton \cite{Wilton} : il y montre que la s\'erie 
\begin{equation}\label{t152}
\psi_1(x)=-\frac{1}{\pi}\sum_{n\soe 1}\frac{\tau(n)}{n}\sin(2\pi nx)\, ,
\end{equation}
o\`u $\tau(n)$ d\'esigne le nombre de diviseurs du nombre entier naturel $n$, converge si, et seulement si la s\'erie $\Wcal(x)$ converge (assertion $(20_{III})$, p. 221). Nous nommerons \emph{points de Wilton} les nombres r\'eels $x$ irrationnels tels que la s\'erie~$\Wcal(x)$ converge, et \emph{fonction de Wilton} la somme de la s\'erie~$\Wcal(x)$ en ses points de convergence. Sa propri\'et\'e fondamentale est l'\'equation fonctionnelle
\begin{equation}
  \label{t93}
  \Wcal (x)=\log (1/x) -x\Wcal(\{1/x\}),
\end{equation}
valable quand $x$ est un point de convergence, $x \neq 0$.

\smallskip

La s\'erie $\psi_1(x)$ n'est autre que celle obtenue formellement \`a partir de la série~$\fhi_1(x)$ en rempla\c{c}ant les termes $B_1(nx)$ par leur d\'e\-ve\-lop\-pe\-ment en s\'erie de Fourier. Le fait que la s\'erie $\fhi_1(x)$ de Chowla converge si, et seulement si la s\'erie $\Wcal(x)$ converge, ne fut d\'emontr\'e qu'en 2004 par la Bret\`eche et Tenenbaum (cf. \cite{breteche-tenenbaum}, Th\'eor\`eme 4.4, p. 16), et, ind\'ependamment, par Oskolkov en 2005 (cf.~\cite{oskolkov}, Corollary 1, A, p. 210). L'identit\'e $\fhi_1(x)=\psi_1(x)$  en tout point de convergence est \'etablie dans ces deux travaux.  
 La Bret\`eche et Tenenbaum d\'eduisent ce dernier r\'esultat de leur th\'eorie g\'en\'erale des s\'eries trigonom\'etriques \`a coefficients arithm\'etiques fond\'ee sur la P-sommation (consistant \`a approcher l'ensemble des nombres entiers par celui des nombres dont tous les diviseurs premiers sont $\ioe y$, avec $y$ tendant vers l'infini), tandis qu'Oskolkov \'etudie une classe g\'en\'erale de proc\'ed\'es de sommation de la s\'erie double sous-jacente \`a~$\fhi_1$ et $\psi_1$.

\smallskip

Notre propos dans le pr\'esent article est, avant tout, m\'ethodologique. Nous allons montrer que l'\'egalit\'e des ensembles de convergence des s\'eries $\fhi_1$ de Chowla et $\Wcal$ de Wilton peut \'egalement \^etre obtenue en suivant la m\'ethode introduite en 1914 par Hardy et Littlewood dans \cite{HL} : celle de l'\'equation fonctionnelle approch\'ee. C'est pr\'ecis\'ement cette m\'ethode que Wilton utilisa en 1933 pour d\'emontrer l'\'egalit\'e des ensembles de convergence de $\psi_1$ et $\Wcal$. 

La m\'ethode de l'\'equation fonctionnelle fournit une expression de la somme $\fhi_1(x)+\demi\Wcal(x)$ en tout point de convergence. L'analyse de cette expression conduit au r\'esultat nouveau suivant.
\begin{theo}
La fonction $G=\fhi_1+\demi\Wcal$ se prolonge en une fonction de p\'eriode~$1$, born\'ee, continue en tout irrationnel, et ayant une discontinuit\'e de premi\`ere esp\`ece en tout rationnel.
\end{theo}

La fonction $G$ est l'objet du \S\ref{t163}, o\`u ses propri\'et\'es sont \'enonc\'ees plus pr\'ecis\'ement. Son caract\`ere born\'e est d\'ej\`a mentionn\'e par Oskolkov dans \cite{oskolkov}, Theorem 1, p. 200. 

Le th\'eor\`eme et sa d\'emonstration ont \'et\'e d\'ecrits dans notre pr\'epublication \cite{autocor} de 2013. Signalons qu'une identit\'e voisine de notre \'enonc\'e, exprim\'ee en termes de la fonction $\psi_1$, a depuis \'et\'e obtenue par Bettin (cf. \cite{bettin}, (1.8), p. 11422).

\smallskip

Nous connaissons deux applications de ce th\'eor\`eme. D'une part, il nous a permis de montrer que l'ensemble des points de d\'erivabilit\'e de la fonction d'autocorr\'elation multiplicative de la fonction {\og partie fractionnaire\fg}, donn\'ee par
\begin{equation}\label{t210}
A(\lambda) = \int_0^{\infty} \{t\} \{\lambda t\} \frac{dt}{t^2} \quad ( \lambda \soe 0),
\end{equation}
et d\'efinie par B\'aez-Duarte \emph{et al.} dans \cite{baez-duarte-all}, est l'ensemble des points de Wilton. Ce r\'esultat, qui sera l'objet d'une publication s\'epar\'ee, a \'et\'e red\'emontr\'e par la Bret\`eche et Tenenbaum dans \cite{bt} \`a partir de r\'esultats obtenus dans \cite{breteche-tenenbaum}.

D'autre part, le th\'eor\`eme et le fait que $\Wcal$ v\'erifie l'\'equation fonctionnelle~\eqref{t93} ont \'et\'e utilis\'es par Maier et Rassias \cite{maier-rassias,maier-rassias2} pour \'etudier le comportement asymptotique des moments de la fonction $\lvert\fhi_1\rvert$.

\smallskip

Voici quel est le plan de cet article. Nous rappelons d'abord au \S\ref{par:fractions-continues} les \'el\'ements de la th\'eorie classique des fractions continues qui nous seront utiles. Le \S\ref{efa} donne un aper\c{c}u historique de la m\'ethode de l'\'equation fonctionnelle approch\'ee. Le c\oe ur de l'article est le \S\ref{t216}, qui pr\'esente un r\'esultat g\'en\'eral de convergence. Le \S\ref{t217} contient l'application de cette th\'eorie aux s\'eries de Chowla et Wilton, et le \S \ref{t218} mentionne succinctement d'autres applications.

\section{Fractions continues}\label{par:fractions-continues} 

\subsection{Notation}

Nous utilisons dans ce paragraphe la notation classique des fractions continues $[b_0; b_1,\dots,b_k]$. \`A partir d'une suite $(b_k)$, on d\'efinit les suites $(p_k)$ et $(q_k)$ par les formules de r\'ecurrence
\begin{align*}
p_{k+1}&=b_{k+1}p_k+p_{k-1}\\
q_{k+1}&=b_{k+1}q_k+q_{k-1} \quad (k \soe 0),
\end{align*}
et les valeurs initiales
\begin{equation}\label{t219}
p_{-1}=1, \, q_{-1}=0 \,; \, p_0 = b_0, \, q_0=1,
\end{equation}
de sorte que
\[
[b_0; b_1,\dots,b_k] = \frac{p_k}{q_k}\cdotp
\]

Dans cet article, $b_0$ sera toujours \'egal \`a $0$.

\subsection{Profondeur d'un nombre réel de $[0,1[$}\label{profondeur_rationnel}

Soit $r$ un nombre rationnel, $0<r<1$, mis sous forme irr\'eductible $p/q$. Il peut s'\'ecrire d'une et une seule fa{\c c}on sous la forme
\begin{equation}\label{u0}
r=[0; b_1,\dots,b_k]
\end{equation}
avec $k \in \Nat^*$, $b_i \in \Nat^*$ pour $1\ioe i\ioe k$, et $b_k \soe 2$. Nous
dirons que $r$ est de profondeur $k$. On a alors $\alpha_k(r)=0$. 

Par convention, $0$ est de profondeur $0$, et tout irrationnel est de profondeur infinie. On voit que, pour tout $x \in\, ]0,1[$ de profondeur $k$, $\alpha(x)$ est de profondeur~$k-1$. 

\subsection{Cellules}\label{par:cellules}

Soit $k\in\Nat^*$, $b_0=0$ et $b_1, \dots, b_k \in\Nat^*$. Par définition, la cellule (de
profondeur~$k$)~$\cgot(b_1,\dots,b_k)$ est l'intervalle ouvert d'extr\'emit\'es
$[b_0;b_1,\dots,b_k]$ et $[b_0;b_1,\dots,b_{k-1},b_k+1]$. Par convention, la seule cellule de profondeur $0$ est $]0,1[$.

Les extr\'emit\'es de $\cgot(b_1,\dots,b_k)$ sont
$$
\frac{p_k}{q_k} \quad \text{et} \quad \frac{p_k+p_{k-1}}{q_k+q_{k-1}}
$$
(dans cet ordre si $k$ est pair ; dans l'ordre oppos\'e si $k$ est impair).

Inversement, soit $r=p/q$ un nombre rationnel (avec $(p,q)=1$), appartenant à $]0,1[$, représenté sous la forme \eqref{u0}. \'Ecrivons $[0; b_1,\dots,b_{k-1}]$ sous forme r\'eduite $p_{k-1}/q_{k-1}$ (si $k=1$, on a~$p_0=0$ et~$q_0=1$). Le nombre rationnel $r$ est une extr\'emit\'e de deux cellules de profondeur $k$ (qui sont donc contig\"ues) :
$$
\cgot \quad \text{d'extr\'emit\'es} \quad [0; b_1,\dots,b_{k-1},b_k-1]=\frac{p-p_{k-1}}{q-q_{k-1}}  \quad \text{et} \quad r \, ;
$$
$$
\cgot' \quad \text{d'extr\'emit\'es} \quad r \quad \text{et} \quad [0; b_1,\dots,,b_{k-1},b_k+1]=\frac{p+p_{k-1}}{q+q_{k-1}}\cdotp
$$

Enfin, observons que les profondeurs des \'el\'ements d'une cellule de profondeur $K$ sont (strictement) sup\'erieures \`a $K$, et que tout nombre de $]0,1[$ de profondeur $>K$ appartient \`a une unique cellule de profondeur $K$.

\subsection{Les fonctions $a_k$, $p_k$, $q_k$} 

Pour tout $k \in \Nat^*$, nous d\'efinissons la fonction \`a valeurs dans $\Nat^*$, $a_k$, sur l'ensemble des nombres $x \in\, ]0,1[$, rationnels ou irrationnels, dont la profondeur est $\soe k$, par la relation
\begin{equation*}
a_k(x) = \lfloor1/\alpha_{k-1}(x)\rfloor. 
\end{equation*}
On pose $a_0(x)=0$ pour tout $ x \in\, ]0,1[$, et on appelle $a_k(x)$ le $k$\up{e} \emph{quotient incomplet} de $x$.

Les fonctions $p_k$, $q_k$ sont d\'efinies pour les m\^emes valeurs de $x$ par le fait que $p_k(x)/q_k(x)$ est la fraction irr\'eductible \'egale \`a $[0;a_1(x),\dots,a_k(x)]$. Les nombres $p_k(x)$ et $q_k(x) $ sont, respectivement, le num\'erateur et le d\'enominateur de la $k$\up{e} fraction r\'eduite de $x$. Les fonctions $p_{-1}$, $q_{-1}$, $p_0$, $q_0$, sont d\'efinies sur $]0,1[$ par \eqref{t219} (avec $b_0=0$).

Dans une cellule de profondeur $k$, $\cgot(b_1,\dots,b_k)$, les fonctions $a_j$, $p_j$, $q_j$ sont donc d\'efinies pour~$j \ioe k+1$, et constantes pour $j \ioe k$ :
$$
a_j(x)=b_j, \quad \frac{p_j(x)}{q_j(x)}=[b_0;b_1,\dots,b_j] \quad (x \in \cgot(b_1,\dots,b_k), \,0\ioe j \ioe k).
$$

\smallskip

Rappelons quels sont les liens entre les fonctions  $p_k$, $q_k$, d'une part, et les fonctions $\alpha_k$, $\beta_k$, d'autre part.
Pour tout $x \in\, ]0,1[$, dont la profondeur est~$\soe k$, on a
\begin{align}
\alpha_k(x)&= -\frac{p_k(x)-xq_k(x)}{p_{k-1}(x)-xq_{k-1}(x)}\label{eq:identite-alpha}  \\
\beta_k(x)&=(-1)^{k-1}\bigl(p_k(x)-xq_k(x)\bigr)=\big|p_k(x)-xq_k(x) \big|. \label{eq:identite-beta}
\end{align}
Si la profondeur de $x$ est $\soe k+1$, on a
\begin{equation}\label{eq:identite-beta2}
 \beta_k(x) =\frac{1}{q_{k+1}(x)+\alpha_{k+1}(x) q_k(x)}\cdotp
\end{equation}
Si $r=p/q \in\, ]0,1[$ est un nombre rationnel de profondeur $K \soe 1$, \'ecrit sous forme irr\'eductible, alors \eqref{eq:identite-beta2} montre que $\beta_{K-1}(r)=1/q$.

En omettant d'indiquer la variable $x$ (de profondeur $\soe k+1$), on d\'eduit de \eqref{eq:identite-beta2} les encadrements 
\begin{equation} \label{encabeta}
 \frac{1}{q_{k+1}+q_k} \ioe \beta_k \ioe \frac{1}{q_{k+1}} \ioe \frac{1}{F_{k+2}}
\quad (k \soe -1),
\end{equation}
o\`u $F_{k+2}$ est le $(k+2)$\up{e} nombre de Fibonacci, et (cf. \cite{Brjuno}, proposition 1),
\begin{equation}\label{enca_gamma}
-\frac{\log (2q_k)}{q_k}\ioe \beta_{k-1}\log (1/\alpha_k)-\frac{\log q_{k+1}}{q_k}\ioe \frac{\log
2}{q_k}\cdotp 
\end{equation}

\smallskip

Pour tout irrationnel $x \in \, ]0,1[$ et tout $s \in \Com$, nous noterons $\rho(x,s)$ le rayon de convergence de la s\'erie enti\`ere
\[
\sum_{k\soe 0} \beta_{k-1}(x)^sz^k.
\]
Si $\sigma=\Re s \soe 0$, on voit que $\rho(x,s) \soe \big((1 +\sqrt{5})/2\big)^{\sigma}$.

\smallskip

Nous utiliserons souvent le fait que, pour tout $x \in\, ]0,1[$, la suite $\big(\beta_k(x)\big)$ est strictement décroissante. Enfin, notons que, si $k \in \Nat$, les fonctions $\alpha_k$ et $\beta_k$ sont continues en tout point de~$]0,1[$ dont la profondeur est $>k$, en particulier en tout point irrationnel.

\section{Aper\c{c}u historique de la m\'ethode de l'\'equation fonctionnelle approch\'ee} \label{efa}

Pour \'etudier l'ordre de grandeur des sommes partielles d'une s\'erie num\'erique dont le terme g\'en\'eral est oscillant, Hardy et 
Littlewood ont d\'evelopp\'e une m\'ethode \'el\'egante et efficace qui consiste \`a exploiter une \'eventuelle \'equation fonctionnelle approch\'ee  v\'erifi\'ee par les sommes partielles de la s\'erie \'etudi\'ee. 
Pour formaliser cette idée, considérons une suite de fonctions $(u_n)$ d\'efinies sur un ensemble $X$ et \`a valeurs complexes ou vectorielles, et notons  \[
 f(x,v) = \sum_{1\le n\le v} u_n(x) \quad(x\in X, v>0). 
 \]  
 ses sommes partielles.
 
 Suivant les valeurs de $x$, on s'interroge sur le comportement asymptotique de $f(x,v)$ lorsque~$v$ tend vers l'infini
 (caract\`ere born\'e, convergence, ordre de grandeur). Il arrive dans certains cas que l'on puisse \'etablir une relation entre les valeurs $f(x,v)$ et $f(T(x), v')$  o\`u $T$ est une certaine application d\'efinie sur $X$ \`a valeurs dans $X$, et $v'$ d\'epend de $v$ et de $x$. Cette relation contient g\'en\'eralement un terme d'erreur d\'ependant de $v$ et $x$, et prend ainsi le nom d'\'equation fonctionnelle approch\'ee. Il s'agit ensuite de l'it\'erer convenablement pour atteindre l'objectif fix\'e. Dans la circonstance favorable o\`u les valeurs successives $v,v',v'',\dots$, r\'esultant de cette it\'eration, d\'ecroissent jusqu'\`a devenir $<1$, on peut remplacer $f(x,v)$ par une somme comportant arbitrairement peu de termes (somme que l'on peut estimer trivialement), voire plus aucun.

\medskip 

L'exemple initialement \'etudi\'e par Hardy et Littlewood en 1914 est celui des sommes partielles de la fonction th\^eta de Jacobi,
\[
\sum_{ n\le v }  \exp(i \pi n^2 x),
\]
et de certaines variantes de cette somme. Il r\'esulte par exemple du th\'eor\`eme 2.128, p. 209 de~\cite{HL} que
\begin{equation}\label{t0}
\sum_{ n\le v }  \exp(i \pi n^2 x)-\frac{e^{i\pi/4}}{\sqrt{x}}\sum_{ n\le xv }  \exp(-i \pi n^2 /x)=O(\sqrt{1/x}) \quad (0<x<1\, , \,v>0).
\end{equation}
De plus, les propri\'et\'es de p\'eriodicit\'e des diverses {\og sommes th\^eta\fg} consid\'er\'ees leur permettent de remplacer la quantit\'e $1/x$ apparaissant dans l'exponentielle complexe de la seconde somme de la relation \eqref{t0} par la partie fractionnaire\footnote{Nous simplifions ici la pr\'esentation de leur raisonnement ; il faut en fait consid\'erer simultan\'ement trois {\og sommes th\^eta\fg} (cf. \cite{HL} pour les d\'etails).}
$\alpha(x)=\{1/x\}$. Avec nos notations, on a donc $T=\alpha$ et $v'=xv$ pour l'exemple \eqref{t0}. L'it\'eration de \eqref{t0} et de ses variantes permet alors \`a Hardy et Littlewood d'obtenir l'estimation
\begin{equation}\label{t1}
\sum_{ n\le v }  \exp(i \pi n^2 x) \ll v\sqrt{\beta_k(x)}+1/\sqrt{\beta_k(x)}, \quad (x\in ]0,1[ \setminus \Rat, \, k\in \Nat).
\end{equation}
On voit que le terme $v\sqrt{\beta_k(x)}$ peut \^etre omis du second membre de \eqref{t1} si l'on a $v\beta_k(x) <1$, ce qui correspond d'ailleurs \`a une somme vide apr\`es $k$ it\'erations. L'estimation \eqref{t1} permet ensuite d'obtenir des informations sur l'ordre de grandeur de la somme partielle \'etudi\'ee en fonction des propri\'et\'es du d\'eveloppement de $x$ en fraction continue, dont on sait qu'il est r\'egi par la suite des~$\alpha_k(x)$. C'est l'un des {\og grands th\'eor\`emes du vingti\`eme si\`ecle\fg} pr\'esent\'es par Choimet et Queff\'elec dans leur livre \cite{choimet-queffelec}.

\smallskip

Plusieurs contemporains de Hardy et Littlewood pr\'ecis\`erent ces r\'esultats ou en \'etudi\`erent diverses g\'en\'eralisations (cf. \cite{vandercorput-1922}, \cite{mordell-1926}, \cite{oppenheim-1928}, \cite{wilton-1927}). En 1933, Wilton \cite{Wilton} mit cette m\'ethode \`a profit pour retrouver certains r\'esultats concernant les sommes partielles 
 \[
\sum_{n\le v } \tau(n) \exp(2i \pi n x) 
\]
\'etudi\'ees un peu plus t\^{o}t par  Chowla \cite{chowla-30,chowla} et Walfisz \cite{walfisz}, et, comme nous l'avons \'evoqu\'e dans l'introduction, donner un crit\`ere de convergence pour les parties r\'eelle et imaginaire de la s\'erie 
\[
\sum_{n\ge 1} \frac{\tau(n)}{n} \exp(2i \pi n x).  
\]  

Dans tous ces exemples, les sommes partielles satisfont \`a une \'equation fonctionnelle approch\'ee o\`u l'application $T$ est la transformation de Gauss. Les majorations ou les crit\`eres de convergence obtenus pour ces s\'eries font donc naturellement intervenir la nature diophantienne de $x$. 

Apr\`es un relatif oubli\footnote{Cf. cependant \cite{chandrasekharan-narasimhan-68} dans le contexte des corps quadratiques.}, cette m\'ethode a connu r\'ecemment un regain d'int\'er\^et, et a \'et\'e notamment utilis\'ee par Rivoal et ses co-auteurs pour traiter une grande vari\'et\'e de situations (cf.~\cite{rivoal-2012}, \cite{rivoal-roques-2013}, \cite{rivoal-seuret-2014}). 

\medskip 

\`A ce stade, il faut signaler que la principale difficult\'e dans la m\'ethode de l'\'equation fonctionnelle approch\'ee consiste pr\'ecis\'ement \`a \'etablir une telle \'equation, qui est souvent le reflet d'une \'equation fonctionnelle exacte satisfaite par une certaine s\'erie g\'en\'eratrice associ\'ee \`a la s\'erie consid\'er\'ee (la fonction th\^eta de Jacobi et la fonction z\^eta de Riemann pour les exemples qui pr\'ec\`edent).  

L'exploitation de cette \'equation fonctionnelle est \textit{a priori} la partie du travail la plus ais\'ee. N\'eanmoins, elle peut parfois s'av\'erer fastidieuse et il serait commode de disposer de certains r\'esultats g\'en\'eraux pour d\'eterminer rapidement les points de convergence de $f(x,v)$. 
Dans ce but, nous proposons au paragraphe suivant l'\'etude d'une famille sp\'ecifique d'\'equations fonctionnelles approch\'ees.

\section{Un crit\`ere de convergence pour les solutions de certaines \'equations fonctionnelles approch\'ees}\label{t216}

\subsection{Forme de l'\'equation fonctionnelle approch\'ee consid\'er\'ee}\label{subsection:hypotheses} 

Soient $a$ un nombre réel strictement positif, $\theta$  et $s$ des nombres complexes. On suppose qu'une application $f: [0,1[ \times ]0,\infty[ \,\vers \Com$ satisfait \`a l'\'equation 
\begin{equation}\label{t200}
f(x,v)-\theta x^s f\big(\alpha(x), x^av \big)= g(x)+ \eps(x,v) \quad (0<x <1, \, v>0),
\end{equation}
et nous \'etudions pour tout $x\in \, [0,1[$, l'existence de 
\[
\Fgot(x)= \lim_{v\to \infty} f(x,v). 
\] 
Une hypothèse assez naturelle pour espérer appliquer avec succès la méthode développée par Hardy et Littlewood est la convergence de $\eps(x,v)$ vers~$0$ lorsque $v$ tend vers l'infini, pour tout $x$ fixé.  

Cette hypothèse est  suffisante pour énoncer un critère quant à l'existence de $\Fgot(x)$ lorsque $x$ est rationnel, mais elle semble insuffisante pour aborder le cas où $x$ est irrationnel. 
Nous proposons ci-dessous un jeu d'hypothèses sous lesquelles un critère d'existence pour $\Fgot(x)$ peut être établi lorsque $x$ est irrationnel.  

\medskip

Nous supposons dans toute la suite que la fonction $\eps$ figurant dans l'équation fonctionnelle approchée \eqref{t200} satisfait l'hypothèse suivante : 

\smallskip

{\bf i)}  pour tout $x \in\, ]0,1[\,$, on a
\[
\eps(x,v) \vers 0 \quad (v \vers \infty). 
\] 

En outre, nous considérons les trois hypothèses suivantes, portant sur les fonctions $f$, $g$ et $\eps$ figurant dans l'équation fonctionnelle approchée~\eqref{t200} :

\smallskip

{\bf ii)} $\;\eps(x,v) \ll 1 \quad(0<x<1, \, x^a v \ge 1)$ ; 

\smallskip

{\bf iii)} $f(x, 1) \ll 1 \quad (0<x<1)$ ; 

\smallskip

{\bf iv)} $f(x, v) \ll 1 +|g(x)| \quad (0<x<1,  \, x^av< 1).$ 

\smallskip


Notons que \eqref{t200} et {\bf i)} entra\^{\i}nent, pour tout $x\in \, ]0,1[$ tel que $\Fgot(x)$ existe, l'existence de $\Fgot(\alpha(x))$ et la relation
\[
\Fgot(x) - \theta x^s \Fgot(\alpha(x)) = g(x). 
\]
C'est l'\emph{\'equation fonctionnelle exacte} associ\'ee \`a \eqref{t200}.

\subsection{Une solution de l'\'equation fonctionnelle exacte associ\'ee}\label{t207}

Pour toute fonction $g : \, ]0,1[ \vers \Com$, on d\'efinit la s\'erie
\begin{equation}\label{t204}
\Scal_{g}(x)= \sum_{j\ge 0} \theta^j \beta_{j-1}(x)^s g\big(\alpha_j(x)\big) \quad (0 < x < 1),
\end{equation}
en convenant que la somme est limit\'ee aux indices $j<K$ si $x$ est un nombre rationnel de profondeur $K$. Nous omettons d'indiquer la d\'ependance en $\theta$ et~$s$, qui sont fix\'es dans toute cette \'etude.

\smallskip

Notons $D_g$ l'ensemble des nombres r\'eels $x\in~]0,1[$ tels que la s\'erie $\Scal_{g}(x)$ converge. Cet ensemble contient donc tous les nombres rationnels de $]0,1[$ et, peut-\^etre, certains nombres irrationnels. Par convention, $0$ n'appartient pas \`a~$D_g$, mais nous prolongeons la d\'efinition de $\Scal_{g}$ en~$x=0$ en posant $\Scal_{g}(0)=0$. 

La s\'erie $\Scal_g$ a \'et\'e d\'efinie de sorte que la proposition suivante soit vraie. 

\begin{prop}\label{t201}
La fonction $S_g$ est solution, sur $D_g$, de l'\'equation fonctionnelle exacte associ\'ee \`a \eqref{t200} :
\[
\Scal_g(x) - \theta x^s \Scal_g\big(\alpha(x)\big) = g(x) \quad ( x \in D_g).
\]
\end{prop}
\dem

Soit $x \in D_g$, et $K$ sa profondeur (\'egale \`a $\infty$, si $x$ est irrationnel). On a
\begin{align*}
\Scal_g(x)- g(x)&=\sum_{0\ioe j < K-1} \theta^{j+1} \beta_{j}(x)^s g\big(\alpha_{j+1}(x)\big)\\
&=\theta x^s \sum_{0\ioe j < K-1} \theta^j \beta_{j-1}\big(\alpha(x)\big)^s g\big(\alpha_j(\alpha(x))\big)\expli{donc $\alpha(x) \in D_g\cup \{0\}$}\\
&= \theta x^s \Scal_g(\alpha(x)),
\end{align*}
car $\alpha(x)$ est de profondeur $K-1$.\fin

\subsection{\'Enoncé du résultat} 

Sous les hypothèses {\bf i)}, {\bf ii)}, {\bf iii)} et {\bf iv)} du \S\ref{subsection:hypotheses}, on a le résultat suivant : 

\begin{prop}\label{prop-general}
 Si $\Fgot(0)=z_0$ existe, alors la limite $\Fgot(x)$ existe pour tout nombre rationnel $x$ de l'intervalle $]0,1]$, et on a  
 \[
 \Fgot(x)=\Scal_{g}(x)+\theta\beta_{K-1}(x)^sz_0, 
 \] 
o\`u $K\in\Nat$ est la profondeur de $x$. Dans le cas contraire, la limite $\Fgot(x)$ n'existe pour aucun nombre rationnel $x$ de $]0,1]$.  

Si $x\in\, ]0,1]$ est irrationnel et tel que $\rho(x,s) > \lvert \theta\rvert$, $\Fgot(x)$ existe si et seulement si la s\'erie $\Scal_{g}(x)$ est convergente, et dans ce cas on a 
 \[
 \Fgot(x)=\Scal_{g}(x).
 \]  
 \end{prop}
Nous d\'emontrons ce r\'esultat dans les sous-paragraphes suivants en distinguant le cas rationnel (ou les hypothèses {\bf ii)}, {\bf iii)} et {\bf iv)} ne seront pas employées) du cas irrationnel. 

\subsection{Le cas où $x$ est rationnel}

Pour ce cas, la seule hypothèse \textbf{i)} suffit. 
Elle entraîne que, pour tout élément $x$ de $]0,1[$, $\Fgot(x)$ existe si et seulement si $\Fgot\big(\alpha(x)\big)$ existe, et dans ce cas, on a 
\begin{equation}\label{eq:wf-exacte}
\Fgot(x)-\theta x^s\Fgot\big(\alpha(x)\big)=g(x). 
\end{equation} 
Nous allons conclure de cette simple remarque que la convergence ou la divergence de $f(0,v)$ conditionne l'existence de $\Fgot(x)$ pour tous les nombres rationnels de $[0,1[$. 

\begin{prop}\label{t215}
S'il existe $z_0\in\Com$ tel que  $\Fgot(0)=z_0$, alors $\Fgot(x)$ existe pour tout nombre rationnel~$x$ de $]0,1]$, et on a 
\begin{equation}\label{eq:formule:rationnel}
\Fgot(x)= \Scal_{g}(x) +\theta^{K} z_0q^{-s},
\end{equation}
o\`u $q$ est le d\'enominateur de la fraction irr\'eductible \'egale \`a $x$, et $K$ la profondeur de $x$. 

Si en revanche, $\Fgot(0)$ n'existe pas, alors il en est de m\^eme pour $\Fgot(x)$ pour tout nombre rationnel~$x$ de $]0,1]$. 
\end{prop}
\dem 

Dans \eqref{eq:formule:rationnel}, nous pouvons \'ecrire $\beta_{K-1}(x)^s$ au lieu de $q^{-s}$.

Lorsque $f(0,v)$ converge, la convergence de $f(x,v)$ et la relation \eqref{eq:formule:rationnel} se d\'emontrent par r\'ecurrence sur la profondeur $K$ de $x$. 

Le cas $K=0$ correspond \`a $x=0$, et d\'ecoule de l'hypoth\`ese faite sur~$\Fgot(0)$.  

Supposons maintenant $K\ge 1$, et le r\'esultat acquis pour tous les nombres rationnels de profondeur $K-1$. Si $x$ est de profondeur $K$, $\alpha(x)$ est de profondeur $K-1$, et on utilise \eqref{eq:wf-exacte} pour en d\'eduire l'existence de $\Fgot(x)$. De plus, on obtient la relation
\begin{align*}
\Fgot(x)&=g(x)+\theta x^s\Fgot\big(\alpha(x)\big)\\
&=g(x) +\theta x^s\Big(\Scal_{g}\big(\alpha(x)\big) +\theta^{K-1} \beta_{K-2}\big(\alpha(x)\big)^sz_0\Big)\expli{hypoth\`ese de r\'ecurrence}\\
&=\Scal_{g}(x) +\theta^{K} \beta_{K-1}(x)^sz_0,
\end{align*}
d'apr\`es la proposition \ref{t201}.

R\'eciproquement, si $x_0$ est un nombre rationnel de profondeur minimale parmi les $x$ tels que~$\Fgot(x)$ existe, la relation \eqref{eq:wf-exacte} entra\^{\i}ne $x_0=0$.\fin

\subsection{Le cas irrationnel} 
 
\subsubsection{It\'eration de l'\'equation fonctionnelle approch\'ee}

Avant de traiter le cas irrationnel, nous explicitons l'\'equation fonctionnelle obtenue par it\'eration de l'identit\'e \eqref{t200}. Le r\'esultat est
\begin{multline}\label{t202}
f(x,v)-\theta^{k} \beta_{k-1}(x)^s f\big(\alpha_k(x), \beta_{k-1}(x)^av \big)=\\
 \sum_{0 \ioe j <k}\theta^{j}\beta_{j-1}(x)^s \Big(g\big(\alpha_j(x)\big)+ \eps\big(\alpha_j(x),\beta_{j-1}(x)^av\big)\Big), 
\end{multline}
\'egalit\'e valable pour tout $x \in [0,1[$ et tout $ k \in \Nat$ tel que $k$ soit inf\'erieur ou \'egal \`a la profondeur de $x$.

Le cas $k=0$ est trivial. Le cas $k=1$ correspond \`a \eqref{t200}.  Maintenant, si~\eqref{t202} est vraie au rang~$k\soe 1$, et si la profondeur de $x$ est sup\'erieure ou \'egale \`a $k+1$, alors l'identit\'e \eqref{t200} appliqu\'ee \`a $\alpha_{k}(x)$ et $ \beta_{k-1}(x)^av$ nous donne 
\begin{equation}\label{t86}
f\big( \alpha_{k}(x), \beta_{k-1}(x)^av\big) - \theta\alpha_{k}(x)^s 
f\big( \alpha_{k+1}(x), \beta_{k}(x)^av \big)
=
g\big(\alpha_{k}(x)\big)+\eps\bigl(\alpha_{k}(x) ,\beta_{k-1}(x)^av\bigr).  
\end{equation}
En ajoutant \eqref{t202} et $\theta^{k}\beta_{k-1}(x)^s\times$\eqref{t86}, on obtient bien la relation \eqref{t202} au rang~$k+1$ au lieu de~$k$.

\smallskip

Pour la solution $\Scal_g$ de l'\'equation fonctionnelle exacte, on d\'emontre de la m\^eme fa\c{c}on que, si~$x \in D_g$ et si le nombre entier $k$ est inf\'erieur ou \'egal \`a la profondeur de $x$, alors $\alpha_k(x) \in D_g\cup\{0\}$ et
\begin{equation}\label{t206}
\Scal_g(x)-\theta^{k} \beta_{k-1}(x)^s \Scal_g\big(\alpha_k(x)\big)=\\
 \sum_{0 \ioe j <k}\theta^{j}\beta_{j-1}(x)^s g\big(\alpha_j(x)\big). 
\end{equation}

\subsubsection{Résultat préliminaire} 
Nous commençons en examinant le comportement de la partie de la somme de \eqref{t202} où intervient la fonction $\eps$. Comme~$x$ est fix\'e dans tout ce paragraphe, nous noterons simplement $\alpha_j$ et~$\beta_j$ pour~$\alpha_j(x)$ et $\beta_j(x)$.

\begin{prop}\label{u214}
On suppose vérifiées les hypothèses {\bf i)}, {\bf ii)}. 
Soit  $x\in\, ]0,1[$ un nombre irrationnel tel que $\rho(x,s) > \lvert \theta\rvert$. On a 
\begin{equation} \label{t91}
 \sum_{0\ioe j <K}\theta^{j}\beta_{j-1}^s\eps(\alpha_j,\beta_{j-1}^a v) \vers 0 \quad (v \vers \infty,\, K\ge 1, \, \beta_{K-1}^a v\ge 1). 
\end{equation} 
\end{prop}
\dem 
Soit $(v_n)_{n\soe 1} $ une suite de nombres r\'eels, et $(K_n)_{n\soe 1} $ une suite de nombres entiers positifs, telles que 
\[
v_n \vers \infty \quad ( n \vers \infty),
\]
et
\[
\beta_{K_n-1}^a v_n \ge 1 \quad(n\ge 1).
\]

On a\footnote{Nous employons ici la notation d'Iverson : $[P]=1$ si la propriété $P$ est satisfaite, et $[P]=0$ sinon.}

\begin{equation*}
 \sum_{0\ioe j <K_n}\theta^{j}\beta_{j-1}^s\eps(\alpha_j,\beta_{j-1}^a v_n) =\sum_{j\soe 0}\theta^{j}\beta_{j-1}^s\eps(\alpha_j,\beta_{j-1}^a v_n)[j < K_n],
\end{equation*} 
avec
\begin{align*}
&\forall j \in \Nat, \quad \theta^{j}\beta_{j-1}^s\eps(\alpha_j,\beta_{j-1}^a v_n)[j < K_n] \vers 0  \quad (n \vers \infty)\expli{d'après {\bf i)}}\\
&\forall j \in \Nat, \quad  \theta^{j}\beta_{j-1}^s\eps(\alpha_j,\beta_{j-1}^a v_n)[j < K_n]  \ll \lvert \theta\rvert^j \beta_{j-1}^{\sigma} \expli{d'après {\bf ii)}}.
\end{align*}
La conclusion \eqref{t91} découle alors du théorème de la convergence dominée de Lebesgue.\fin

\subsubsection{Condition nécessaire pour la convergence de $f(x,v)$} 

\begin{prop}\label{prop:CN}
On suppose vérifiées les hypothèses {\bf i)}, {\bf ii)} et {\bf iii)}. 
 Soit $x\in\, ]0,1[$ un nombre irrationnel tel que $\rho(x,s) > \lvert \theta\rvert$. Si la  limite $\Fgot(x)$ existe alors la s\'erie $\Scal_{g}(x)$ est convergente et on a 
 \[
 \Fgot(x)= \Scal_{g}(x). 
 \] 
\end{prop}
\dem
Pour $K \in \Nat$, on pose \[v = v(x,K)=1/ (\beta_{K-1}^a v)\] de sorte que
$v(x,K) \vers \infty$ quand $K\vers \infty$. En appliquant la relation \eqref{t202} à~$k = K$, on obtient 
\begin{equation}
\label{t111}
f(x,v) - \theta^K\beta_{K-1}^s  f(\alpha_{K},1) 
=  \sum_{0 \ioe j <K}\theta^{j}\beta_{j-1}^s g\big(\alpha_j\big) 
+ \sum_{0\ioe j <K}\theta^{j}\beta_{j-1}^s\eps(\alpha_j,\beta_{j-1}^a v). 
\end{equation}

D'après {\bf iii)}, 
\begin{equation}\label{t92}
\theta^K\beta_{K-1}^s  f(\alpha_{K},1)  \ll \lvert  \theta\rvert^K\beta_{K-1}^{\sigma}. 
\end{equation} 
 Comme $\lvert \theta \rvert < \rho(x,s)$, le majorant de \eqref{t92} tend vers $0$ quand $K$ tend vers l'infini.
Si $f(x,v)$ converge vers $\Fgot(x)$, les relations \eqref{t91}, \eqref{t111} et~\eqref{t92} m\`enent \`a la conclusion que la s\'erie $\Scal_{g}(x)$
converge et a pour somme $\Fgot(x)$.
\fin

\subsubsection{Condition suffisante pour la convergence de $f(x,v)$}

\begin{prop}\label{prop:CS}
On suppose vérifiées les hypothèses {\bf i)}, {\bf ii)} et {\bf iv)}. 
 Soit $x\in\, ]0,1[$ un nombre irrationnel tel que $\rho(x,s) > \lvert \theta\rvert$. Si la  s\'erie $\Scal_{g}(x)$ est convergente alors la limite $\Fgot(x)$ existe  et on~a 
 \[
 \Fgot(x)= \Scal_{g}(x). 
 \] 
  \end{prop}
\dem 
Pour $v>1$, d\'efinissons l'entier~$K=~K(x,v)\soe~0$ par l'encadrement
$$
\beta_{K-1}^a v \soe 1 > \beta_{K}^a v.
$$
On voit que $K(x,v)\vers \infty$ si, et seulement si $v\vers \infty$.

L'équation \eqref{t111} reste valable. La condition ${\bf iv)}$ entraîne que 
\begin{equation}\label{t203}
\theta^K\beta_{K-1}^s  f(\alpha_{K},\beta_{K-1}^a v)   \ll \lvert  \theta\rvert^K\beta_{K-1}^{\sigma}\big(1+\lvert g\big(\alpha_{K}\big)\rvert \big).
\end{equation} 
Comme $\lvert \theta\rvert < \rho(x,s)$, les relations \eqref{t91}, \eqref{t111}, \eqref{t203} et la convergence de $\Scal_{g}(x)$  entra\^{\i}nent la convergence de $f(x,v)$ vers~$\Scal_{g}(x)$ lorsque $v$ tend vers l'infini.\fin

\subsubsection{Conclusion} 

La proposition \ref{prop-general} résulte des propositions \ref{t215}, \ref{prop:CN} et \ref{prop:CS}.

\subsection{\'Etude de $\Scal_{g}$ quand $g$ est continue sur $[0,1]$}\label{t211}

Dans ce paragraphe, nous supposons que $\theta$ et $s$ v\'erifient
\begin{equation}\label{t205}
\sigma=\Re s > 0 \quad ; \quad \lvert \theta \rvert < \big((1+\sqrt{5})/2\big)^{\sigma}.
\end{equation}

Si $g$ est born\'ee sur $]0,1[$, la s\'erie d\'efinie par \eqref{t204} converge alors normalement sur $]0,1[$, o\`u sa somme $\Scal_g$ est donc d\'efinie et born\'ee. Rappelons que nous avons prolong\'e $\Scal_g$ en $0$ en posant~$\Scal_g(0)=0$.

\begin{prop}\label{prop:gcontinue}
On suppose que la fonction $g$ se prolonge contin\^{u}ment \`a l'intervalle $[0,1]$. 
La fonction~$\Scal_g : \, ]0,1[ \, \to \Com$ est alors continue en chaque irrationnel. De plus, 
\begin{align*}
\Scal_g(x) & \vers g(0) & (x\vers 0, \, x>0)\;\\
\Scal_g(x) & \vers g(1)+ \theta g(0) & (x\vers 1, \, x <1),
\end{align*}
et, pour tout nombre rationnel $r\in\, ]0,1[$ de profondeur $K\ge1$ et dont le d\'enominateur de la fraction irr\'eductible est $q$, les limites \`a droite et \`a gauche de $\Scal_g(x)$ quand $x$ tend vers $r$ sont
\[
\Scal_g(r)+q^{-s}\theta^{K} g(0) \quad  \text{ et } \quad  \Scal_g(r)+q^{-s}\theta^{K}\big(g(1)+\theta g(0)\big),
\]
dans cet ordre si $K$ est pair, dans l'ordre inverse, si $K$ est impair.
\end{prop}
\dem 

Pour $j \in \Nat$ et $0\ioe x<1$, posons
\[
\mu_j(x)=
\begin{cases}
\theta^j \beta_{j-1}(x)^s g\big(\alpha_j(x)\big) & \text{si $x$ est de profondeur $>j$},\\
\qquad\qquad 0 & \text{si $x$ est de profondeur $\ioe j$}.
\end{cases}
\]
On a donc
\begin{equation}\label{t220}
\Scal_g(x)=\sum_{j \soe 0} \mu_j(x) \quad (0\ioe x <1).
\end{equation}

Comme $g$ est born\'ee, il d\'ecoule de \eqref{t205} que la s\'erie \eqref{t220} est normalement convergente sur~$[0,1[$. De plus, la fonction $\mu_j$ est continue en tout point de profondeur $>j$, en particulier en tout irrationnel. La fonction $\Scal_g$ est donc continue en tout irrationnel.

\medskip 

On a ensuite
\begin{align*}
\mu_j(x) &\vers
\begin{cases}
g(0) & (j=0)\\
0 & (j>0)
\end{cases}
&(x \vers 0, \, x>0)\\
\mu_j(x) &\vers
\begin{cases}
g(1) & (j=0)\\
\theta g(0) & (j=1)\\
0 & (j>1)
\end{cases}
& (x \vers 1, \, x <1)
\end{align*}
donc 
\begin{align*}
\Scal_g(x) &\vers g(0) &(x \vers 0, \, x>0)\;\\
\Scal_g(x) &\vers g(1) +\theta g(0) &(x \vers 1, \, x<1).
\end{align*}

\smallskip

\`A pr\'esent, soit    
\[
r=\frac pq=[0;a_1,\ldots,a_K] \quad ((p,q)=1, \, a_K\soe 2).
\] 
un nombre rationnel de profondeur $K\ge 1$. Consid\'erons les cellules contig\"ues de profondeurs
respectives  $K$ et  $K+1$, auxquelles $r$ est adh\'erent :
\[
 \cgot=\cgot(a_1,\ldots,a_K) \quad ; \quad
\cgot'=\cgot(a_1,\ldots,a_K-1,1).
\]
Si $K$ est pair (resp. 
impair), alors la cellule
$\cgot$ se situe \`a droite (resp. \`a gauche) de $r$. 
Nous allons montrer que 
\begin{equation*} 
\lim_{ \substack{x\to r \\ x\in\cgot} } \Scal_g(x)=\Scal_g(r)+\theta^{K}g(0)q^{-s}
\end{equation*}
tandis que 
\begin{equation*}
\lim_{ \substack{x\to r \\ x\in\cgot'} } \Scal_g(x)=\Scal_g(r)+\theta^{K}\big( g(1)+\theta g(0) \big)q^{-s}.
\end{equation*} 

Remarquons tout d'abord que tout nombre rationnel appartenant \`a la r\'eunion 
$\cgot\cup \cgot'$ est de profondeur sup\'erieure ou \'egale \`a $K+1$. 
Pour $j< K$, les fonctions $\alpha_j$ et $\beta_j$ sont continues sur~$\overline{\cgot\cup \cgot'}$.  

\smallskip

$\bullet$ Calcul de $\lim_{ \substack{x\to r \\ x\in\cgot} } \Scal_g(x)$

Pour $x\in\cgot$, on a $p_K(x)=p$, $q_K(x)=q$, et par cons\'equent, d'apr\`es \eqref{eq:identite-alpha}, 
\begin{align*}
\alpha_K(x)&=-\frac{p-x q}{p_{K-1}-x q_{K-1}}=o(1) & (x\to r, \, x\in\cgot)\\
\beta_K(x)&=\lvert p-x q\rvert=o(1) & (x\to r, \, x\in\cgot).
\end{align*}

D'apr\`es \eqref{t206}, nous avons 
\begin{align*}
 \Scal_g(x)=
\sum_{j\le K} \theta^j\beta_{j-1}(x)^s g\big(\alpha_j(x)\big)
+\theta^{K+1}\beta_{K}(x)^s\Scal_g\big(\alpha_{K+1}(x)\big)\quad(x\in\cgot).
\end{align*}

Comme $\Scal_g$ est born\'ee sur $[0,1[$, il suit  
\begin{equation}\label{eq:limite-cgot}
\lim_{ \substack{x\to r \\ x\in\cgot} } \Scal_g(x)
=\sum_{j< K} \theta^j \beta_{j-1}(r)^s g\big(\alpha_j(r)\big)
+\theta^{K}\beta_{K-1}(r)^sg(0)= 
\Scal_g(r)+\theta^{K}g(0)q^{-s}.
\end{equation}

$\bullet$ Calcul de $\lim_{ \substack{x\to r \\ x\in\cgot'} } \Scal_g(x)$

Tous les nombres rationnels de $\cgot'$ ont une profondeur sup\'erieure ou \'egale \`a 
$K+2$. Pour $x\in\cgot'$, on a $p_K(x)=p-p_{K-1}$, $q_K(x)=q-q_{K-1}$, $p_{K+1}(x)=p$ et $q_{K+1}(x)=q$, de sorte que 
\begin{align*}
\alpha_K(x)&=-\frac{p-p_{K-1}-x(q-q_{K-1})}{p_{K-1}-x q_{K-1}}=1+o(1) &( x\to r, \, x\in\cgot')\\
\alpha_{K+1}(x)&=-\frac{p-x q}{p-p_{K-1}-x(q-q_{K-1})}=o(1) &(x\to r, \, x\in\cgot')\\
 \beta_K(x)&=\lvert p-p_{K-1}-x(q-q_{K-1})\rvert= \frac 1q+o(1) & (x\to r, \, x\in\cgot') .
\end{align*}

D'apr\`es \eqref{t206}, nous avons 
\begin{align*}
 \Scal_g(x)=
\sum_{j\le K} \theta^j\beta_{j-1}(x)^s g\big(\alpha_j(x)\big)
+\theta^{K+1}\beta_{K}(x)^s\Scal_g\big(\alpha_{K+1}(x)\big)\quad(x\in\cgot').
\end{align*}

Il suit 
\begin{align*} 
\lim_{\substack{x\to r\\ x\in\cgot'}}  \Scal_g(x)&= 
 \sum_{j< K}  \theta^j\beta_{j-1}(r)^s g\big(\alpha_j(r)\big)+
\theta^{K}q^{-s}g(1)+\theta^{K+1}q^{-s}g(0)\\
&=\Scal_g(r)+\theta^{K}\big(g(1)+\theta g(0)\big)q^{-s}.   \fine
\end{align*}

\section{Application aux s\'eries de Chowla et de Wilton}\label{t217}

\subsection{La fonction de Wilton}\label{t164}

La s\'erie de Wilton \eqref{t209} n'est autre que \eqref{t204} dans le cas particulier
\[
\theta=-1 \quad  ; \quad s=1 \quad ; \quad g=\log \, :
\]
\[
\Wcal(x)=\sum_{k\ge 0}(-1)^k \beta_{k-1}(x)\log\big(1/\alpha_k(x)\big),
\]
avec la convention du \S\ref{t207} pour les nombres $x$ rationnels. Sa somme, la fonction de Wilton v\'erifie donc l'\'equation fonctionnelle \eqref{t93} sur l'ensemble de convergence, constitu\'e des nombres rationnels et des nombres de Wilton. Le fait qu'un nombre irrationnel $x$ soit un nombre de Wilton peut s'exprimer en termes des d\'enominateurs des r\'eduites de $x$ ; il \'equivaut \`a la convergence de la s\'erie
\begin{equation*}
  \sum_{k\ge 0} (-1)^k \frac{\log q_{k+1}(x)}{q_k(x)}\cdotp
\end{equation*}
C'est une cons\'equence directe de l'encadrement \eqref{enca_gamma}
et du fait que la s\'erie $$\sum_{k\ge 0} \frac{\log q_k}{q_k}$$ est convergente.

\smallskip

Si l'on prend $\theta=1$ au lieu de $\theta=-1$ ci-dessus, on obtient la s\'erie et la fonction de Brjuno,
\begin{equation}
\label{def-brjuno}
\Phi(x)=\sum_{k\ge 0} \beta_{k-1}(x)\log\big(1/\alpha_k(x)\big),
\end{equation}
solution de l'\'equation fonctionnelle
\begin{equation}
  \label{t96}
  \Phi (x)=\log (1/x) +x\Phi(\{1/x\}),
\end{equation}
qui ne diff\`ere de \eqref{t93} que par un signe. 
La fonction de Brjuno, qui joue un r\^ole important dans la th\'eorie de certains syst\`emes dynamiques (cf. \cite{yoccoz-1995}), appara\^{\i}t d\'ej\`a en filigrane dans l'article de Wilton (cf. \cite{Wilton} (7.32), p. 235).

\subsection{La fonction  $G$}\label{t163}

Posons pour $x >0$,
\begin{equation}\label{eq:def-F} 
 F(x)=\frac{x+1}{2}A(1)-A(x)-\frac{x}{2} \log x,
\end{equation} 
o\`u $A$ est la fonction d'autocorr\'elation d\'efinie par \eqref{t210}. La fonction $F$ est continue sur $]0,\infty[$, s'annule en $x=1$, et se prolonge par continuit\'e en
$0$ en posant\footnote{On a $A(1)=\log 2\pi-\gamma$, o\`u $\gamma$ d\'esigne la constante d'Euler (cf. d\'emonstration de la proposition 8 de \cite{baez-duarte-all}).} $F(0)=A(1)/2$. 

Consid\'erons ensuite la s\'erie $\Scal_F$ et sa somme. Les r\'esultats du \S\ref{t211} s'appliquent : la fonction~$\Scal_F$ est d\'efinie et born\'ee sur $[0,1[$, continue en tout irrationnel, et ayant des limites \`a droite et \`a gauche en tout rationnel, comme \'enonc\'e dans la proposition \ref{prop:gcontinue}. On d\'efinit ensuite la fonction $G$ par p\'eriodicit\'e \`a partir de $\Scal_F$ :
\[
G(x) =\Scal_F(\{x\}) \quad (x \in \Real).
\]

En r\'esum\'e :

$\bullet$ la fonction $G$, de p\'eriode $1$, est d\'efinie par
\begin{equation}\label{eq:serieF}
 G(x)= \sum_{j\ge 0}(-1)^j \beta_{j-1}(x)
  F\big(\alpha_j(x)\big) \quad (0 <x <1),
\end{equation}
avec la convention du \S\ref{t207} pour les nombres $x$ rationnels, en particulier on a~$G(0)=0$ ;

$\bullet$ la fonction $G$ v\'erifie l'\'equation fonctionnelle 

  \begin{equation}
    \label{t98}
G(x)=F(x)-x G\bigl (\alpha(x)\bigr)\quad (0 <x <1) \, ;   
  \end{equation}

$\bullet$  si $r=p/q\in\, ]0,1[$ est un nombre rationnel de profondeur $K\soe 1$, \'ecrit sous forme irr\'eductible,  alors $G$ a des limites \`a droite et \`a gauche au point $r$, donn\'ees par
\[
G(r)+(-1)^{K} \frac{A(1)}{2q} \quad  \text{ et } \quad  G(r)+(-1)^{K+1} \frac{A(1)}{2q}\virg
\]
dans cet ordre si $K$ est pair, dans l'ordre inverse, si $K$ est impair. Pour $K=0$, on obtient les limites \`a droite en $0$ et \`a gauche en $1$.

En particulier, $G$ est \emph{normalis\'ee} : sa valeur en tout point est la moyenne de ses limites \`a gauche et \`a droite.

\subsection{\'Equation fonctionnelle approch\'ee pour les sommes partielles de~$\fhi_1(x)$}\label{par:fhi1}

Les sommes partielles de $\fhi_1$ v\'erifient une \'equation fonctionnelle {\og a\-symp\-to\-ti\-que\-ment proche\fg} de celle v\'erifi\'ee par $G-\Wcal/2$.

\begin{prop}\label{prop:eq-approchee-fhi1}
Pour $0<x<1$ et $v >0$, on a
$$
\sum_{1\ioe m\ioe v} \frac{ B_1(m x)}{m} +x \sum_{1\ioe n\ioe x
v} \frac{ B_1\bigl (n\alpha(x)\bigr)}{n}=F(x)-\demi\log
(1/x)+\eps_1(x,v) ,$$
avec 
\begin{equation}\label{t213}
\eps_1(x,v)\ll \frac{1}{xv} \quad (xv\soe 1). 
\end{equation} 
\end{prop}
\dem Nous reproduisons la formule figurant en haut de la page 225 de \cite{baez-duarte-all} :
\begin{multline}\label{t119}
\sum_{1\ioe m\ioe v} \frac{ B_1(m x)}{m} +x \sum_{1\ioe n\ioe x
v} \frac{ B_1 (n/x)}{n}=\\
\frac{x}{2}\int_0^v\{t\}^2t^{-2}dt +
\frac{1}{2}\int_0^{x v}\{t\}^2t^{-2}dt- \int_0^v\{t\}\{x
t\}t^{-2}dt
+\frac{x -1}{2} \log (1/x)\\+\frac{x -1}{2}\int_v^{x v}\{t\}t^{-2}dt+\frac{1}{2 x v}(\{x v\}-x\{v\})^2 +\frac{x -1}{2x v}(\{x v\}-x\{v\}).
\end{multline}
Comme $B_1 (n/x)= B_1\bigl (n\alpha(x)\bigr)$, on reconna\^it bien la formule annonc\'ee, avec
\begin{multline*}
  \eps_1(x,v)=-\frac{x}{2}\int_v^{\infty}\{t\}^2t^{-2}dt-\demi\int_{x v}^{\infty}\{t\}^2t^{-2}dt+\int_v^{\infty}\{t\}\{x t\}t^{-2}dt\\
+\frac{x -1}{2}\int_v^{x v}\{t\}t^{-2}dt+\frac{1}{2 x v}(\{x v\}-x\{v\})^2 +\frac{x -1}{2 x v}(\{x v\}-x\{v\}).
\end{multline*}
Chacun des six termes composant $\eps_1(x,v)$ est $\ll 1/x v$ si $v\soe 1/x$ (ce qui entra\^ine $v\soe 1$).\fin

\medskip

Signalons que la source de l'identit\'e \eqref{t119} est la troisi\`eme d\'emonstration de la loi de r\'eciprocit\'e quadratique, propos\'ee par Gauss en 1808 (cf. \cite{gauss}, \S 5). Eisenstein en 1844 (cf. \cite{028.0828cj}) en donna une pr\'esentation g\'eom\'etrique particuli\`erement intuitive \`a l'aide d'un comptage de points \`a coordonn\'ees enti\`eres dans un rectangle. L'article d'Eisenstein fut traduit par Cayley  et publi\'e en 1857 dans le \textit{Quarterly Journal of pure and applied Mathematics}, \'edit\'e par Sylvester (cf. \cite{cayley}). Trois ans plus tard, Sylvester publia une note au Comptes Rendus \cite{sylvester} o\`u il g\'en\'eralisait l'argument d'Eisenstein \`a un rectangle~$[0,v]\times [0,x v]$ o\`u $v$ et $x$ sont des r\'eels quelconques. La relation \eqref{t119} d\'ecoule de l'identit\'e de Sylvester par sommation partielle (cf. \cite{baez-duarte-all}, p. 223-225).

\subsection{Conclusion}\label{t161}

Posons
\[
f(x,v)=-2\sum_{1\ioe m\ioe v} \frac{ B_1(m x)}{m} +2G(x) \quad (0 \ioe x <1, \, v>0).
\]
La proposition \ref{prop:eq-approchee-fhi1} et la relation \eqref{t98} fournissent l'\'equation fonctionnelle approch\'ee
\[
f(x,v) + x f\big(\alpha(x), xv \big)= \log(1/x)+ \eps(x,v) \quad (0<x <1, \, v>0),
\]
avec $\eps=-2\eps_1$. 
On est donc dans le cadre du \S\ref{subsection:hypotheses}, avec 
\[
\theta=-1\; ; \; s=1\; ;\; a=1\; ;\; g(x)=\log(1/x). 
\]
On déduit immédiatement de la majoration $\eps(x,v) \ll (xv)^{-1}$ (pour $xv\ge 1$), que les hypoth\`eses~{\bf i)} et {\bf ii)}  sont v\'erifi\'ees.  L'hypoth\`ese {\bf iii)} r\'esulte du fait que~$G$ et $B_1$ sont born\'ees, et  l'estimation~$\sum_{n \ioe v}1/n \ll 1+\log v$ (pour~$v \soe 1$) entraîne {\bf iv)}. La proposition \ref{prop-general} s'applique donc. On en conclut que les s\'eries~$\fhi_1(x)$ et $\W(x)$ convergent pour les m\^emes valeurs de $x\in[0,1[$ et  que l'identit\'e
\begin{equation}\label{eq:identite-G-X}
\fhi_1(x) = -\frac12 \W(x) + G(x) 
\end{equation} 
est valable en tout point de convergence.  Compte tenu des propri\'et\'es de $G$ vues au \S\ref{t163}, cela termine la d\'emonstration du th\'eor\`eme.

\section{Autres applications}\label{t218}

Dans cette section, nous montrons comment les résultats du 
\S\ref{t216}
 permettent de retrouver et parfois de préciser quelques résultats antérieurs pour certaines séries oscillantes. 

\subsection{Les séries  $\sum_{n\ge 1} \frac{\tau(n)}{n}\sin(2\pi nx)$
et $\sum_{n\ge 1} \frac{\tau(n)}{n}\cos(2\pi nx)$
} 
Dans \cite{Wilton}, Wilton établit l'existence d'une fonction $\mathfrak{F}:[0,1] \to \Com$ continue, et telle que pour~$0<~x~\ioe~ 1$ et $x^2v\soe 1$, on a
\begin{equation}\label{EFA-psi1}
\sum_{n\ioe v} \frac{\tau(n)}{n}\sin(2\pi nx)
+x\sum_{n\ioe x^2v} \frac{\tau(n)}{n}\sin (2\pi n/x) =
 \frac{\pi}2\log(1/x) +\Im(\mathfrak{F}(x)) +O\big((x^2v)^{-1/5}\big)
\end{equation}
et
\begin{multline}\label{EFA-psi2}
\sum_{n\ioe v} \frac{\tau(n)}{n}\cos(2\pi nx)-x\sum_{n\ioe x^2v} 
\frac{\tau(n)}{n}\cos(2\pi n/x) =\\
\demi\log^2(1/x) 
+(\gamma-\log 2\pi)\log (1/x) 
+\Re(\mathfrak{F}(x)) +O\big((x^2v)^{-1/5}\big).
\end{multline}
Dans l'article \cite{balazard-martin-2015}, où nous avons élaboré une nouvelle démonstration des relations \eqref{EFA-psi1} et \eqref{EFA-psi2},  nous avons obtenu une expression de la fonction $\mathfrak{F}$ en fonction de $A$. En particulier, on a pour tout $x\in[0,1]$, 
\begin{equation}
\label{eq:id-FF}
\Im(\mathfrak{F}(x)) = -\pi F(x).  
\end{equation}

Comme évoqué\ dans l'introduction, Wilton déduit de \eqref{EFA-psi1} que les points de convergence de la série $\psi_1(x)$ définie en \eqref{t152} sont les nombres rationnels et les nombres de Wilton. Les éléments réunis au \S\ref{t216} permettent de retrouver ce résultat et également d'obtenir une nouvelle démonstration de l'identité~$\fhi_1(x) =\psi_1(x)$. Posons 
\[
h(x,v) = \frac2\pi \sum_{n\ioe v} \frac{\tau(n)}{n}\sin(2\pi nx) +2 G(x), 
\] 
où la fonction $G$ est définie au \S\ref{t163}. 
D'après \eqref{t98}, \eqref{EFA-psi1} et \eqref{eq:id-FF}, on a\[
h(x,v)+xh(\alpha(x),x^2v) = \log(1/x) + \eps(x,v) \quad(0<x<1, v>0)
\]  
avec $\eps(x,v)= O((x^2v)^{-1/5})$.
Nous sommes dans le cadre du \S\ref{subsection:hypotheses}, avec 
\[
\theta=-1\; ; \; s=1\; ;\; a=2\; ;\; g(x)=\log(1/x). 
\]
Les hypoth\`eses {\bf i)}, {\bf ii)} et {\bf iii)} sont, comme au \S\ref{t161}, immédiatement vérifiées.  Quant à l'hypoth\`ese~{\bf iv)} elle r\'esulte du fait que $G$ est born\'ee et de l'estimation
uniforme de Walfisz (cf.~\cite{walfisz},~$(25_{VI})$, p.~566),  
\begin{equation*}
\sum_{n\ioe v} \frac{\tau(n)}{n}\sin(2\pi nx)=
\sum_{m\le v} \frac{1}{m}\sum_{k\ioe v/m}\frac{1}{k} \sin(2\pi km x)\ll 1+ \log (v) \quad(v\ge 1, \, x\in\, ]0,1[).
 \end{equation*} 
La proposition \ref{prop-general} s'applique et montre que les 
s\'eries $\psi_1(x)$ et $\W(x)$ convergent pour les m\^emes valeurs de $x\in[0,1[$ et que l'identit\'e
\begin{equation}\label{eq:identite-G-psi}
\psi_1(x) = -\frac12 \W(x) + G(x) 
\end{equation} 
est valable en tout point de convergence. Compte tenu de \eqref{eq:identite-G-X}, nous obtenons bien $\psi_1(x)= \fhi_1(x)$ en tout nombre rationnel et tout nombre de Wilton.   

Un raisonnement similaire permet de retrouver le résultat de Wilton selon lequel pour tout $x\in\,]0,1[$, la série 
\[
\psi_2(x) = \sum_{n\ge 1} \frac{\tau(n)}{n}\cos(2\pi nx)
\] 
converge si et seulement si $x$ est irrationnel et satisfait à 
\begin{equation}\label{serie-log-carre}
\sum_{k\ge 0} \beta_{k-1}(x) \log^2(1/\alpha_k(x))<\infty. 
\end{equation}
On peut montrer que cette dernière condition est équivalente à 
\[
\sum_{k\ge 0} \frac{\log^2(q_{k+1}(x))}{q_k(x)}<\infty. 
\] 
De plus, si l'on pose pour tout $x$ irrationnel de $]0,1[$ 
\[
\Phi_2(x) = \sum_{k\ge 0} \beta_{k-1}(x) v(\alpha_k(x)), 
\] 
avec $v(x) = \frac{1}{2}\log^2(1/x)+ (\gamma-\log(2\pi))\log(1/x)$, les résultats du \S\ref{t216}  permettent de déduire de la continuité de $\Re(\mathfrak{F})$ que la fonction 
$
\psi_2 - \Phi_2
$, définie en tout nombre irrationnel $x\in\, ]0,1[$ tel que~\eqref{serie-log-carre}, est bornée,  
et peut être prolongée  à une fonction définie sur $[0,1]$, continue en tout nombre irrationnel et admettant une limite à droite à gauche en tout nombre rationnel.

\subsection{Une s\'erie \'etudi\'ee par Rivoal et Roques}
La série 
\[
\Psi(x) = \sum_{n\ge 1}\frac{\sin(2\pi n^2 x)\cot(\pi nx)}{n^2} 
\]
 a été introduite et étudiée par Rivoal et Roques dans l'article \cite{rivoal-roques-2013}. Ils remarquent que $\Psi(x)$ est divergente en tout nombre $x$ rationnel et obtiennent (theorem 2, p. 100 de \cite{rivoal-roques-2013})  que pour tout nombre irrationnel $x\in\, ]0,1[$, la série $\Psi(x)$ est convergente si et seulement si $x$ est un nombre de Brjuno, c'est-à-dire est tel que $\Phi(x)<\infty$ où $\Phi(x)$ est la fonction de Brjuno définie par \eqref{def-brjuno}.  
Compte tenu de \eqref{enca_gamma}, c'est encore équivalent à 
\begin{equation}
\label{equiv-RR}
\sum_{k\ge 0} \frac{\log q_{k+1}(x)}{q_k(x)}<\infty. 
\end{equation}

Rivoal et Roques déduisent cette caractérisation d'une équation fonctionnelle approchée pour les sommes partielles 
\[
\Psi(x,v)=\sum_{1\le n\ioe v} \frac{\sin(2\pi n^2 x)\cot(\pi nx)}{n^2}\virg
\]
qui prend la forme suivante (theorem 1 de \cite{rivoal-roques-2013}) :  pour tous $x\in\, ]0,1]$, $v > 0$, on a 
\begin{equation}
\label{EFA-RR}
\Psi(x,v) - x\Psi
(1/x,xv) = \frac{1}{\pi x}
\sum_{1\le m\le v} \frac{\sin(2\pi m^2 x)}{m^3}
+ G(x,v),  
\end{equation}
où $G:\,]0,1] \times ]0,\infty[\to \Real$ est telle que : 
\begin{itemize}
\item pour tout $x\in\,]0,1]$, $\lim_{v\to \infty} G(x,v)$ existe, on la note $G(x)$  ; 
\item il existe $C>0$ tel que pour tout $(x,v) \in \,]0,1] \times ]0,\infty[$, $|G(x,v)|\le C$. 
\end{itemize}
\smallskip 

Des estimations élémentaires conduisent à la formule
\[
\frac{1}{\pi x}\sum_{1\le m\le v} \frac{\sin(2\pi m^2 x)}{m^3} 
= \log(1/x) + H(x) + O\big((xv^2)^{-1}\big), 
\]
valable pour $x>0$ et $v >0$ et 
où $H : \, ]0,\infty[\to \Real$ est bornée. La série $\mathcal{S}_{H+G}(x)$ avec $\theta =s=1$ (cf. \eqref{t204}) définit une fonction bornée sur $]0,1[$, et 
\[
F(x,v) = \Psi(x,v) - \mathcal{S}_{H+G}(x)
\] 
satisfait à l'équation fonctionnelle approchée 
\[
F(x,v) - xF(1/x,xv) = \log(1/x)+\eps(x,v)  
\] 
avec $\eps(x,v)  =G(x,v)-G(x) + O\big((xv^2)^{-1}\big)$. 
Nous sommes donc dans le cadre du \S\ref{subsection:hypotheses} avec 
\[
\theta=1\quad ; \quad  s=1\quad ; \quad  a=1 \quad \text{et} \quad g(x)=\log(1/x). 
\]
L'expression de~$\eps(x,v)$ montre que {\bf i)} et {\bf  ii)} sont satisfaites. L'hypothèse {\bf iii)} est vérifiée car la quantité $\sin(2\pi x) \cot(\pi x)=2\cos^2 (\pi x)$ est bornée.  Quant à l'hypothèse~{\bf iv)}, elle découle de la majoration 
\[
\norm{nt}\le n\norm{t} \quad(t\in\Real,\, n\in\Nat)
\]
où $\norm{x}$ désigne l'entier le plus proche du nombre réel $x$. 
Comme $\sin(\pi x) \asymp~\norm{x}$, on en déduit en effet que pour $0<x<1$, $xv <1$,  on a  
\[
F(x,v) \ll 1+ \sum_{1\le n\le v } \frac{1}{n^2} \frac{\norm{n^2 x}}{\norm{nx}} 
 \le  1+\sum_{1\le n\le 1/x } \frac{1}{n} \ll 1+ g(x). 
\]  
Notons enfin que la série $\Psi(0)$ est divergente. 
La proposition \ref{prop-general} redonne bien le critère obtenu par Rivoal et Roques pour la convergence de la série $\Psi(x)$, et le fait que la fonction 
$\Psi - \mathcal{S}_{H+G}$, définie en tout nombre de Brjuno, est bornée.

\begin{center}
  {\sc Remerciements}
\end{center}
\begin{quote}
{\footnotesize Outre leurs laboratoires respectifs, les auteurs remercient les UMI 2615 (CNRS-UIM, Moscou) et~3457 (CNRS-CRM, Montréal), qui ont fourni des conditions idéales pour leur travail sur cet article.
}
\end{quote}

\medskip

\footnotesize

\noindent BALAZARD, Michel\\
Aix Marseille Univ, CNRS, Centrale Marseille, I2M, Marseille, France\\
Adresse \'electronique : \texttt{balazard@math.cnrs.fr}

\medskip

\noindent MARTIN, Bruno\\
ULCO,  LMPA, Calais, France\\
Adresse \'electronique : \texttt{Bruno.Martin@univ-littoral.fr}

\end{document}